\documentclass[12pt,bezier,epsf]{article}
\usepackage{amsmath,amssymb,amsfonts}
\usepackage{float,latexsym,shortvrb}
\usepackage{graphicx}
\usepackage[usenames]{color}
\usepackage{ulem}
\numberwithin{equation}{section}
\def\cU{{\cal U}}
\def\cT{{\cal T}}
\def\cV{{\cal V}}
\def\cS{{\cal S}}
\def\cB{{\cal B}}
\def\cW{{\cal W}}
\def\cD{{\cal D}}
\def\sse{\subseteq}
\def\lora{\longrightarrow}
\def\lan{\langle}
\def\ran{\rangle}
\def\la{\leftarrow}
\def\Ra{\Rightarrow}
\def\ra{\rightarrow}

\def\ds{\displaystyle}
\def\rul{\rule{2mm}{2mm}}

\title{ON $\cU$-EQUIVALENCE SPACES}
\author{{ Farshad  Omidi$^{a}$, MohammadReza Molaei$^{b}$} \\ {\tiny $^{a}$Department of Mathematics, Velayat  University, Iranshahr, Iran   }\\
{\tiny E-mail: faromidimath66@gmail.com}\\ {\tiny $^{b}$Department of Pure Mathematics, Shahid Bahonar University of Kerman , Kerman, Iran}\\ {\tiny and Center of Excellence on Modeling and Control Systems}\\ {\tiny Ferdowsi University of Mashhad, Mashhad, Iran}\\
{\tiny E-mail: mrmolaei@uk.ac.ir}}
\date{}

\begin{document}
\maketitle

\begin{abstract}
In this paper induced $\cU$-equivalence spaces are introduced and discussed. Also the notion of $\cU$-equivalently open subsets of a $\cU$-equivalence space and $\cU$-equivalently open functions are studied. Finally, equivalently uniformisable topological spaces are considered.
\end{abstract}
\section{Introduction}

Non-Archimedean spaces have been considered from different viewpoints \cite{B, BR, D, DI, G, MO, P, R1, R2}.
In this paper we  deal with ${\cU}$-equivalence spaces. These spaces have been introduced first in 2014  \cite{OM}. A ${\cU}$-equivalence space $(X,{\cU})$ is a set $X$ along with a collection ${\cU}$ of equivalence relations on $X$ such that ${\cU}$ is closed under finite intersections. A {\cU}-equivalence space is a structure close to uniform spaces \cite{J, K, M} and fuzzy uniform spaces \cite{H, L}.\\ A function $f:X\lora Y$ where $(X,{\cU})$ and $(Y,{\cV})$ are two ${\cU}$-equivalence spaces is called ${\cU}$-equivalently continuous if $(f\times f)^{-1}(V)\in{\cU}$ whenever $V\in{\cV}$, where
\[(f\times f)^{-1}(V)=\{(x,y)\in X\times X | (f(x),f(y))\in V\}.\]
Moreover if $(X,{\cU})$ is a $\cU$-equivalence space, then the collection $\cT_{\cU}=\{G\sse X| $ for each $x\in G$, there exists $U\in\cU$ such that $U[x]\sse  G\}$ is a topology on $X$ with the base $\{U[x] | U\in\cU, x\in X\}$ where $U[x]=\{y\in X | (x,y)\in U\}$. We refer to $\cT_{\cU}$ as the $\cU$-induced topology. We are going to consider induced $\cU$- equivalence classes and $\cU$-products. Also we will introduce and discuss $\cU$-equivalently open subspaces,  and equivalently uniformisable spaces.

\section{Induced $\cU$-equivalence spaces and $\cU$-products}

Let $\{\Phi_i : X\lora X_i\}_{i\in I}$ be an indexed family of  functions where $X$ is a set and for each $i\in I$, \ \  $(X_i,\cU_i)$ is a $\cU$-equivalence space. The idea is to induce a $\cU$-equivalence class \cite{OM} on $X$ for which each $\Phi_i$ is $\cU$-equivalently continuous without making the $\cU$-equivalence class on $X$ unnecessarily strong.

{\bf Definition 2.1.}
Let $\cS$ be a family of equivalence relations on a set $X$. Then the collection of all finite intersections of members of $\cS$ (that forms a $\cU$-equivalence class on $X$) called the $\cU$-equivalence class generated by $\cS$ and it is denoted by $\lan\cS\ran$.

Note that $\lan\cS\ran$ is the smallest $\cU$-equivalence class on $X$ which contains $\cS$.

{\bf Proposition 2.2.}
Let $\{\varphi_i: X\lora X_i | i\in I\}$ be an indexed family of functions where $X$ is a set and for each $i\in I$, $(X_i,\cU_i)$ is a $\cU$-equivalence space. Then there  exists a smallest $\cU$-equivalence class on $X$ for which each $\Phi_i$ is $\cU$-equivalently continuous.

{\bf Proof:}
Let $\cS^{\la}=\{(\Phi_i\times\Phi_i)^{-1}(U_i) | U_i\in \cU_i, i\in I\}$.
  Obviously, for each $U_i\in\cU_i$, $i\in I$, \ \ $(\Phi_i\times\Phi_i)^{-1}(U_i)$ is an equivalence relation on $X$. Let $\cU^{\la}$ be the $\cU$-equivalence class generated by $\cS^{\la}$. Clearly for each $i\in I$, $\Phi_i$ is $\cU$-equivalently continuous w.r.t. this $\cU$-equivalence class.

Finally, if $\cU$ is a $\cU$-equivalence class on $X$ for which each $\Phi_i$ is $\cU$-equivalently continuous and $U\in \cU^{\la}$, then $U=\ds\bigcap^n_{j=1}(\Phi_{i_j}\times\Phi_{i_j})^{-1}(U_{i_j}), U_{i_j}\in \cU_{i_j}$. Since each $\Phi_i$ is $\cU$-equivalently continuous w.r.t. $\cU$ and $\cU$ is closed under finite intersections, then $U\in\cU$. Hence $\cU^{\la}\sse\cU$.~$\rul$

The $\cU$-equivalence class $\cU^{\la}$ in the last proposition is called the induced $\cU$-equivalence class. Note that if each $\cU_i$ has a base $\cB_i$, then the collection $\{(\Phi_i\times\Phi_i)^{-1}(U_i) | U_i\in \cB_i, i\in I\}$ generates $\cU^{\la}$.

{\bf   Corollary 2.3.}
Let $\Phi:X\lora Y$ be a function where $X$ is a set and $(Y,\cV)$ is a $\cU$-equivalence space.

Then $\cV^{\la}= \{(\Phi\times\Phi)^{-1}(V) | V\in \cV\}$.

The following property is a characteristic of the induced $\cU$- equivalence classes.

{\bf Proposition 2.4.}
Let $\Phi:X\lora Y$ and $\Psi:Y\lora Z$ be functions where  $(X,\cU)$, $(Y,Z)$ and $(Z,\cW)$ are $\cU$- equivalence spaces. If $Y$ has the induced $\cU$-equivalence class, then $\Phi$ is $\cU$-equivalently continuous if and only if $\Psi\Phi$ is $\cU$-equivalently continuous.

{\bf Proof:}
The proposition follows from the equality
\[(\Psi o\Phi\times\Psi o\Phi)^{-1}(W)=(\Phi\times\Phi)^{-1}((\Psi\times\Psi)^{-1}(W))\]
where $W$ runs through the members of $\cW$ and the fact that $\cV=\cW^{\la}$.~$\rul$

Let $(X,\cU)$ be a $\cU$-equivalence space and let $A\sse X$. By Corollary 2.3. (let $\Phi$ to be inclusion map) the collection $\{A\times A\cap U|U\in\cU\}$ is the induced $\cU$-equivalence class on $A$ that is called relative $\cU$-equivalence class and denoted by $\cU/A$ (see \cite{OM}).

Hence $\cU/A$ is the smallest $\cU$-equivalence class on $A$ makes inclusion map $\cU$-equivalently continuous.

{\bf Definition 2.5.}
(a)   Let $(X,\cU)$ and $(Y,\cV)$ be two $\cU$-equivalence spaces. A bijection $\Phi:X\lora Y$ is said to be a $\cU$-equivalence (function) if $\Phi$ and $\Phi^{-1}$ are $\cU$-equivalently continuous.
(b)  A function $f:X\lora Y$ is said to be $\cU$-embedding if it is one to one and a $\cU$-equivalence when regarded as a function  from $(X,\cU)$ on to $(f(X),\cV/f(X))$.

{\bf Theorem 2.6.}
Let $(X,\cU)$ and $(Y,\cV)$ be two $\cU$-equivalence spaces and let \linebreak $\Phi:X\lora Y$ be a function. Then the following statements are equivalent:
(a)   $\Phi$ is a $\cU$-embedding.
(b)  $\Phi$ is one to one, $\cU$-equivalently continuous and $\cU=\cV^{\la}$

{\bf Proof:}
(a)$\Ra$(b). Since $\Phi$ is $\cU$-equivalently continuous when regarded as a function from $(X,\cU)$, then $\cV^{\la}\subseteq\cU$.
  Conversely, if $U\in\cU$, then $(\Phi^{-1}\times\Phi^{-1})^{-1}(U)\in\cV / \Phi(X)$.
  Hence $(\Phi^{-1}\times\Phi^{-1})^{-1}(U)=(\Phi(X)\times\Phi(X))\cap V$ for some $V\in\cV$.
So $U=(\Phi\times\Phi)^{-1}(V)$. Hence $\cU\sse\cV^{\la}$.

(b)$\Ra$(a). Since $\cU=\cV^{\la}$, then $\Phi$ is $\cU$-equivalently continuous. To complete the proof, we will prove that $\Phi^{-1}:(\Phi(X),\cV/\Phi(X))\lora(X,\cU)$ is $\cU$-equivalently continuous. For this, let $U\in\cU$, $V\in\cV$ and $U=(\Phi\times\Phi)^{-1}(V)$.

So $(\Phi^{-1}\times\Phi^{-1})^{-1}(U)=(\Phi(X)\times\Phi(X))\cap V$. Hence $\Phi^{-1}$ is $\cU$-equivalently continuous.~$\rul$

Not all $\cU$-equivalently continuous  injections are $\cU$-embedding's. For example let $X$ be a set and $\cU$ be the collection of all equivalence relations on $X$ and $\cV$ consists of $X^2$ and let $\Phi$ be the identity map on $X$. Then $\Phi$ is a $\cU$-equivalently injection  when regarded as a function from $(X,\cU)$ on to $(X,\cV)$ but not $\cU$-embedding if $|X|>1$.
A useful sufficient condition is

{\bf  Proposition 2.7.}
Let $\Phi:X\lora Y$ be a $\cU$-equivalently continuous function where $(X,\cU)$ and $(Y,\cV)$ are $\cU$-equivalence spaces. Suppose $\Phi$ admits a $\cU$-equivalently left inverse. Then $\Phi$ is a $\cU$-embedding.

{\bf Proof:}
Suppose $\Psi:Y\lora X$ is the $\cU$-equivalently left inverse of $\Phi$. By theorem 2.6 it is sufficient to show $\cU=\cV^{\la}$. Suppose $U\in\cU$ and $V=(\Psi\times\Psi)^{-1}(U)$. Since  $\Psi$ is $\cU$-equivalently continuous, then $V\in\cV$ and $(\Phi\times\Phi)^{-1}(V)=(\Phi\times\Phi)^{-1}((\Psi\times\Psi)^{-1}(U))=(\Psi o\Phi\times\Psi o\Phi)^{-1}(U)=U$. The last equality is true because $\Psi o\Phi$ is the identity map on $X$. Hence $\cU=\cV^{\la}$.$\rul$

{\bf  Proposition 2.8.}
Let $\{X_j\}$ be a finite covering of the $\cU$-equivalence space $(X,\cU)$.

Suppose that $X_j$ is totally bounded \cite{OM} for each index $j$. Then $X$ is totally bounded.

{\bf  Proof:}
We recall that $X$ is totally bounded if for each $U\in\cU$ there exists $x_1,x_2,\dots,x_n\in X$ such that $X=\ds\bigcup^n_{i=1} U[x_i]$. Let $X=\ds\bigcup^m_{j=1} X_j$ and let $U\in\cU$ and $U_j=X_j\times X_j\cap U$ for each $j$.

Since $X_j$ is totally bounded, then there exists $x_1^j,x_2^j,\dots,x_{n_j}^j\in X_j$ such that \linebreak $X_j=\ds\bigcup^{n_j}_{i=1}U_j[x_i^j]$.
Hence $X=\ds\bigcup^{m}_{j=1}\ds\bigcup^{n_j}_{i=1}U[x_i^j]$.$\rul$

The converse of Proposition 2.8 is true. In fact we have

{\bf Proposition 2.9.}
Let $\Phi:X\lora Y$ be a function where $X$ is a set and $(Y,\cV)$ is a $\cU$-equivalence space.
If $Y$ is totally bounded, then so is $X$ with the induced $\cU$-equivalence class.

{\bf Proof:}
Let $V\in\cV$ and let $U=(\Phi\times\Phi)^{-1}(V)$. There exist $y_1, y_2,\dots,y_n\in Y$ such that $Y=\ds\bigcup^n_{i=1}V[y_i]$. If $T=\{y_i|\Phi^{-1}(V[y_i])\neq\Phi\}$, then $T\neq\Phi$. For each $y_i\in T$, choose a member of the non-empty set $\Phi^{-1}(V[y_i])$, say, $x_i$. Suppose $x\in X$ and let $y=\Phi(x)$.
Then $y\in V[y_i]$ for some $i$. Hence $(y_i,y)\in V$. Thus $(\Phi(x_i),y)\in V$ because $V$ is an equivalence relation on $X$. Consequently $X=\ds\bigcup^m_{i=1}U[x_i]$.$\rul$

The notation of induced $\cU$-equivalence classes is the standard way of giving a $\cU$-equivalence class to the cartesian product $X=\Pi X_i$ where $\{(X_i,\cU_i) | i\in I\}$ is  an indexed family of $\cU$-equivalence spaces.

More precisely, let $\{(X_i,\cU_i) | i\in\}$ be given.
For each $i\in I$ let $\pi_i$ be the $i$-th canonical projection.
By Proposition 2.2, there exists a smallest $\cU$-equivalence class $\Pi\cU_i$ for which $\pi_i$ is $\cU$-equivalently continuous. The $\cU$-equivalence space $(\Pi X_i,\Pi\cU_i)$ is called the $\cU$-product of family $\{(X_i,\cU_i) | i\in I\}$.

{\bf  Proposition 2.10.}
Let $\{(X_i,\cU_i) | i\in I\}$ be a family of $\cU$-equivalence spaces and let \linebreak $\Phi:A\lora\Pi X_i$ be a function where $A$ is a $\cU$-equivalence space and $(\Pi X_i,\Pi\cU_i)$ is the \linebreak $\cU$-product. Then $\Phi$ is $\cU$-equivalently continuous if and only if each of the functions \linebreak $\Phi_j=\pi_j o\Phi:A\lora X_j$ is $\cU$-equivalently continuous.

{\bf  Proof:}
The proof is an easy consequence of the equality $(\pi_j o\Phi\times \pi_j o\Phi)^{-1}(U_j)$ \linebreak $=(\Phi\times\Phi)^{-1} ((\pi_j\times\pi_j)^{-1}(U_j))$ for $U_j\in\cU_j$, $j\in I$.$\rul$

According to \cite{OM} a $\cU$-equivalence space $(X,\cU)$ is separated if the intersection of all members of $\cU$ coincides with $\Delta X=\{(x, x) | x\in X\}$.

{\bf Proposition 2.11.}
Let $\{(X_i,\cU_i) | i\in I\}$ be a family of separated  $\cU$-equivalence spaces.
Then the $\cU$-product is separated .

{\bf  Proof:}
Let $x=(x_i)$ and $y=(y_i)$ be two members of $X=\Pi X_i$ and $(x,y)\in U$ for all $U\in\Pi\cU_i$.

So $(x,y)\in(\pi_i\times\pi_i)^{-1}(U_i)$ for each $U_i\in\cU_i$, $i\in I$. Hence $x=y$. Consequently the intersection of all members of $\Pi\cU_i$ coincides with $\Delta X$. Hence $(\Pi X_i,\Pi \cU_i)$ is separated.$\rul$

{\bf Theorem 2.12.}
Let $\{(X_i,\cU_i) | i\in I\}$ be a family of separated  $\cU$-equivalence spaces.
Then $\cT_{_{\Pi\cU_i}}=\Pi\cT_{_{\cU_i}}$.

{\bf  Proof:}
Let $j\in I$ be fixed and let $X=\Pi X_i$, We contend $\pi_j:X\lora X_j$ is continuous when regarded as a function from topological space $(X,\cT_{_{\Pi\cU_i}})$ onto topological space $(X_j,\cT_{_{\cU_j}})$. To see this, let $G_j\in\cT_{_{\cU_j}}$ and $G=\pi^{-1}_j(G_j)$. Suppose $x=(x_i)\in G$. Thus $x_j\in G_j$. So, there exists $U_j\in \cU_j$ such that $U_j[x_j]\sse G_j$.

 If $U=(\pi_j\times\pi_j)^{-1}(U_j)$, then $U\in\Pi\cU_i$. Suppose $z=(z_i)\in U[x]$. Thus $z_j\in U_j[x_j]\sse G_j$. So $z\in G$.
  Thus $U[x]\sse G$. This proves $\pi_j$ is continuous. Hence by definition of $\Pi\cT_{_{\cU_i}}$, $\Pi\cT_{_{\cU_i}}\sse \cT_{_{\Pi\cU_i}}$.
  For other way inclusion, suppose $G\in\cT_{_{\Pi\cU_i}}$ and $x\in G$. Then there exists $U_{i_j}\in\cU_{i_j}$ for $1\leq j\leq n$ such that $\ds\bigcap^n_{j=1}(\pi_{i_j}\times\pi_{i_j})^{-1}(U_{i_j})\sse U$ and $U[x]\sse G$ for some $U\in\Pi\cU_i$.

 Assume $G_{i_j}=U_{i_j}[x_{i_j}]$ for $1\leq j\leq n$. Then $G_{i_j}\in\cT_{_{\cU_{i_j}}}$ and $x\in\ds\bigcap^n_{j=1}\pi_{i_j}^{-1}(G_{i_j})\sse G$. Hence $\cT_{_{\Pi\cU_i}}\sse \Pi\cT_{_{\cU_i}}$. Consequently,  $\cT_{_{\Pi\cU_i}}= \Pi\cT_{_{\cU_i}}$.

 \section{$\cU$-Equivalently Open Subspaces}

 Among the subspaces of a $\cU$-equivalence space $(X,\cU)$ special attention should be given to those which are $\cU$-equivalently open, in the sense that the inclusion map is $\cU$-equivalently open. We recall that a function $f:X\lora Y$ where $(X,\cU)$ and $(Y,\cV)$ are $\cU$-equivalence spaces is called $\cU$-equivalently open if for each $U\in\cU$, there exists $V\in\cV$ such that $V[f(x)]\sse f(U[x])$ for all $x\in X$. Roughly speaking, every $\cU$-equivalently open function is locally surjective (see \cite{OM}). For example $\Phi$ is always $\cU$-equivalently open when $\cV$ is discrete i.e. $\cV$ is the collection of all equivalence relations on $Y$.

  Let $(X,\cU)$ be a $\cU$-equivalence space. Then $\cU$ is called rich if $X^2\in\cU$.

{\bf  Proposition 3.1.}
Let $(X,\cU)$ be a rich $\cU$-equivalence space and let $A\sse X$. Then the following statements are equivalent:

(a) $A$ is $\cU$-equivalently open.

(b) for each $U\in\cU$, there exists $V\in\cU$ such that $V[x]\sse U[x]\cap A$ for all points $x\in A$.

(c) there exists $V_0\in\cU$ such that $V_0[x]\sse A$ for all points $x\in A$ (in particular $A\in\cT_{_{\cU}}$).

{\bf  Proof:}
(a)$\ra$(b). Let $U\in\cU$ and $\tilde{U}=A\times A\cap U$.

Then $\tilde{U}\in\cU/A$. So there exists $V\in\cU$ such that $V[x]\sse\tilde{U}[x]$ for all points $x\in A$.

Hence $V[x]\sse U[x]\cap A$ for all $x\in A$.

(b)$\ra$(c). Since $X^2\in\cU$, then there exists $V_0\in\cU$ such that $V_0[x]\sse X^2[x]\cap A$ for all $x\in A$. Hence $V_0[x]\sse A$ for all points $x\in A$.

(c)$\ra$(a). Let $\tilde{U}=A\times A\cap U\in\cU|A$ and let  $V=V_0\cap U$. Then $V\in\cU$ and $V[i(x)]\sse i(\tilde{U}[x])$ for all points $x\in A$. Hence $A$ is $\cU$-equivalently open.$\rul$

Suppose $\alpha,\beta:X\lora Y$ are maps from $X$ into $Y$.
The coincidence set of $\alpha$ and $\beta$ is the set $C(\alpha,\beta)=\{x\in X | \alpha(x)=\beta(x)\}$.
Also if $\varphi:X\lora Y$ is a function from the $\cU$-equivalence space $(X,\cU)$ into the set $Y$, let us say that $\varphi$ is transverse to $X$ (or $(X,\cU)$) if $(\varphi\times\varphi)^{-1}(\Delta Y)\cap U=\Delta X$ for some $U\in\cU$.
Roughly speaking, $\varphi$ is transverse to $X$ if $\varphi$ is one to one on a region of $X^2$.

{\bf Proposition 3.2.}
Let $(X,\cU)$ , $(Y,\cV)$ and $(Z,\cW)$ be $\cU$-equivalence spaces, \linebreak $\alpha,\beta:X\lora Y$ be two $\cU$-equivalently continuous functions and $\varphi:Y\lora Z$ be transverse to $Y$. If $\varphi o\alpha=\varphi o\beta$ and $\cU$ is rich then the coincidence set $C(\alpha,\beta)$ is  $\cU$-equivalently open in $X$.

{\bf Proof:}
To prove the proposition, we use part (c) of proposition (3.1). By the transitivity of $\varphi$, there exists $V\in\cV$ such that $(\varphi\times\varphi)^{-1}(\Delta Z)\cap V=\Delta Y$.

It is easy to see that $U_0[x]\sse C(\alpha,\beta)$ for all $x\in C(\alpha,\beta)$ where \[U_0=(\alpha\times\alpha)^{-1}(V)\cap(\beta\times\beta)^{-1}(V).~\rul\]

{\bf Definition 3.3.}
Let $f:X\lora Y$ be a function where $X$ is a set and $(Y,\cV)$ is a $\cU$-equivalence space. $f$ is called $\cU$-equivalently surjection if for each $y\in Y$ and $V\in\cV$, there exists $x\in X$ such that $y\in V[f(x)]$.

Obviously any surjection from $X$ onto the $\cU$-equivalence space $(Y,\cV)$ is a $\cU$-equivalently surjection.

On the other hand if $X=Y=\bf{R}$ (the set of real numbers) and $\cV=\{\bf{R}^2\}$, we define $f:X\lora Y$ by $f(x)=x^2$, then $f$ is a $\cU$-equivalently surjection  but not a surjection.

{\bf Proposition 3.4.}
Let $f:X\lora Y$ be a $\cU$-equivalently surjection. If $f$ is $\cU$-equivalently open and $\cU$ is rich, then $f$ is a surjection.

{\bf Proof:}
Since $f$ is $\cU$-equivalently open, then there exists $V\in\cV$ such that $V[f(x)]\sse f(X)$ for all $x\in X$ because $\cU$ is rich. If $y\in Y$, then there exists $x\in X$ such that $y\in V[f(x)]$.
Hence $y\in f(X)$ that means $f$ is a surjection.$\rul$

{\bf Definition 3.5.}
Let $(X,\cU)$ be a $\cU$-equivalence space. A subset $D$ of $X$ is said to be $\cU$-equivalently dense (in $X$), if for each $x\in X$ and each $U\in\cU$, there exists $a\in D$ such that $(a,x)\in U$.

Note that $D$ is $\cU$-equivalently dense iff $D$ is a dense subset of topological space $(X,\cT_{_{\cU}})$.
Also, note that $D$ is $\cU$-equivalently dense in $X$ iff the inclusion map $i:D\ra X$ is a $\cU$-equivalently surjection.

{\bf Proposition 3.6.}
Let $(X,\cU)$ be a rich $\cU$-equivalence space and let $A$ be a  $\cU$-equivalently open and $\cU$-equivalently dense subset of $X$. then $A$ equals $X$.

{\bf Proof:}
Since $A$ is $\cU$-equivalently open, then the inclusion map $i:A\lora X$ is $\cU$-equivalently open. If $x\in X$ and $U\in\cU$, then there exists $a\in A$ such that $(a,x)\in U$. So $x\in U[i(a)]$. Hence the inclusion map is also a $\cU$-equivalently surjection. Now Proposition 3.4 implies the inclusion map is a surjection. Hence $A$ equals $X$.$\rul$

{\bf Corollary 3.7.}
Let $(X,\cU)$, $(Y,\cV)$ and $(Z,\cW)$ be $\cU$-equivalence spaces and let \linebreak $\alpha,\beta:X\lora Y$ be $\cU$-equivalently continuous and $\varphi:Y\lora Z$ be transverse to $Y$. If $\varphi o\alpha=\varphi o\beta$, $\cU$ is rich and if the coincidence set $C(\alpha,\beta)$ is $\cU$-equivalently dense, then $\alpha=\beta$.

We say the $\cU$-equivalence space $(X,\cU)$ is connected if the topological space $(X,\cT_{_{\cU}})$ is connected. We have the following strange result:

{\bf  Proposition 3.8.} Every non-empty subset of a connected $\cU$-equivalence space is $\cU$-equivalently dense.

{\bf Proof:}
Let $A$ be a non-empty subset of a connected  $\cU$-equivalence space $(X,\cU)$ and let $U\in\cU$. Since  $A$ is non-empty, then so is $U[A]$. We contend that $U[A]$ is clopen. To see this, since $U[A]=\ds\bigcup_{x\in A} U[x]$, then $U[A]$ is open. On the other hand, $\overline{U[A]}=\ds\bigcup_{v\in\cU} V[U[A]]$, where $\overline{U[A]}$ is the closure of $U[A]$ w.r.t. $\cT_{\cU}$. Hence
\[\overline{U[A]}\sse U[U[A]]\sse U[A].\]
The last statement is true because $U$ is an equivalence relation on $X$. So $U[A]$ is closed.  Hence the conectedness of the $\cU$-equivalence space $(X,\cU)$ implies that $U[A]=X$. Now if $x\in X$, then $x\in U[A]$.
Thus $(a,x)\in U$ for some $a\in A$ that means $A$ is $\cU$-equivalently dense in $X$.$\rul$

{\bf Corollary 3.9.} Every $\cU$-equivalently open subset of a connected  $\cU$-equivalence space is either empty or full.

{\bf Proof:} immediate from (3.6) and (3.8).$\rul$

{\bf Corollary 3.10.} Let $(X,\cU)$, $(Y,\cV)$ and $(Z,\cW)$ be $\cU$-equivalence spaces and let \linebreak $\alpha,\beta:X\lora Y$ be $\cU$-equivalently continuous and $\varphi : Y\lora Z$  be transverse to $Y$. If  $\varphi o\alpha=\varphi o\beta$, $\cU$ is rich and $\alpha(x_0)=\beta(x_0)$ for some point $x_0$ in $X$ and if $(X,\cU)$ is connected, then $\alpha=\beta$.

{\bf Proof:}   The coincidence set $C(\alpha,\beta)$ is non-empty because $\alpha(x_0)=\beta(x_0)$. Now the result follows from Corollary 3.7 and 3.9.$\rul$

\section{Equivalently uniformisable spaces}

In this section we consider topological spaces which are equivalently uniformisable. More precisely, let $(X,\cT)$ be a topological space. The question is:  under what condition(s) there exists a $\cU$-equivalence class $\cU$ on $X$ such that $\cT_{_{\cU}}=\cT$. We begin with the following definition.

{\bf Definition 4.1.}
A topological space $(X,\cT)$ is said to be equivalently uniformisable if there exists a $\cU$-equivalence class on $X$ such that $\cT=\cT_{_{\cU}}$.

{\bf Proposition 4.2.}
Let $(X,\cT)$ and $(Y,\cV)$ be respectively topological and $\cU$-equivalence  spaces and let $f:X\lora Y$ be a topological equivalence from topological  space $(X,\cT)$  on to topological space $(Y,\cT_{_{\cV}})$.
Then the  topological space $(X,\cT)$   is equivalently uniformisable.

{\bf Proof:}
Let $\cU$ be the induced $\cU$-equivalence class  on $X$ by $\cV$ (i.e. $\cU=\overleftarrow{\cV}=\{(f\times f)^{-1}(V) | V\in\cV\}$. We contend that $\cT=\cT_{_{\cU}}$.
Suppose $G\in\cT$ and $x\in G$. Since $f$ is a topological equivalence, then $f(G)\in \cT_{_{\cV}}$. Hence there exists $V\in\cV$ such that $V[f(x)]\sse f(G)$. Let $U=(f\times f)^{-1}(V)$. Then $U\in\cU$ and $U[x]\sse G$. Hence $\cT\sse \cT_{_{\cU}}$.
Conversely let $G\in \cT_{_{\cU}}$ and  $x\in G$. Since $x\in G$ and $G\in\cT_{_{\cU}}$, then there exists $U\in\cU$ and $V\in\cV$ with $U=(f\times f)^{-1}(V)$ and $U[x]\sse G$. Let $O=V[f(x)]$ and $G'=f^{-1}(O)$. Then $G'\in\cT$ and $x\in G'\sse G$. This shows that $\cT_{_{\cU}}\sse\cT$. Consequently $\cT=\cT_{\cU}$. This means that topological space $(X,\cT)$ is equivalently uniformisable.$\rul$

{\bf Proposition 4.3.}
Let $(X,\cU)$ be a  $\cU$-equivalence  space and let $Y\sse X$. Then $\cT_{_{\cU/Y}}=\cT_{_{\cU}}/Y$.

{\bf Proof:}
We recall that $\cT_{_{\cU}}/Y=\{Y\cap G | G\in \cT_{_{\cU}}\}$ and
\[\cT_{_{\cU/Y}}=\{G\sse Y | \forall y\in G, \exists V_y\in\cU / Y, V_y[y]\sse G\}\]
and $\cU/Y = \{Y\times Y\cap U | U\in \cU\}$. We now turn to the proof.
We first show that $\cT_{_{\cU/Y}}\sse\cT_{_{\cU}}/Y$.  Let $G\in \cT_{_{\cU/Y}}$ and let $y\in G$. Then there exists $V_y\in \cU/Y$ and $U_y\in\cU$ such that $V_y=Y\times Y\cap U_y$ and $V_y[y]\sse G$.
Let $G_1=\ds\bigcup_{y\in G} U_y[y]$. Then $G_1\in\cT_{_{\cU}}$ and $G=Y\cap G_1$. So $\cT_{_{\cU/Y}}\sse\cT_{_{\cU}}/Y$.
Now let $G\in \cT_{_{\cU}}/Y$ and let $y\in G$. Then $G=Y\cap G_1$ for some $G_1\in\cT_{_{\cU}}$. Since $G_1\in \cT_{_{\cU}}$, then there exists $U_y\in\cU$ such that $U_y[y]\sse G_1$. If $V_y=Y\times Y\cap U_y$, then $V_y[y]\sse G$. This shows that $\cT_{_{\cU}}/Y\sse \cT_{_{\cU/Y}}$. Hence these two topologies are the same.$\rul$

{\bf Proposition 4.4.}
Let $(X,\cT)$ and $(Y,\cV)$ be respectively topological and  $\cU$-equivalence  spaces and let $f:X\lora Y$ be a topological embedding from $(X,\cT)$ into topological space $(Y,\cT_{_{\cV}})$. Then $(X,\cT)$ is equivalently uniformisable.

{\bf Proof.} The function $f:X\lora Z$ is a topological equivalence where $Z=f(X)$ regarded as a topological space with $\cT_{_{\cV}}/Z$ (the induced topology on $Z$ by $\cT_{_{\cV}}$). But by Proposition 4.3, $\cT_{_{\cV}}/Z=\cT_{_{\cV/Z}}$. Hence by Proposition 4.2, $(X,\cT)$ is equivalently uniformisable.$\rul$

By Proposition 4.4, a topological space is equivalently uniformisable if it can be embedded into a topological space whose topology induced by a $\cU$-equivalence class.

{\bf Definition 4.5.} Let $(X,d)$ be a pseudo metric space and let $r$ be a positive real number. Then

(a) $d$ is called $r$-transitive if $d(x,y)<r$, $d(y,z)<r$ implies $d(x,z)<r$ for all points $x$, $y$ and $z$ in $X$.

(b) $d$ is called transitive if for each $r>0$, $d$ is $r$-transitive.

For an example of $r$-transitive pseudo-metric, let $X$ be a non-empty set and let
\[d_\alpha(x,y)=\left\{\begin{array}{cc} ~~~~\alpha~~~~ & x\neq y \\ 0 & x=y \end{array}\right.\]
where $\alpha$ is a positive real number. Then for each $r>\alpha$, $d_\alpha$ is $r$-transitive.

Let $d$ be a transitive pseudo-metric  on $a$ set $X$ and let $$B_d(r)=\{(x,y)\in X\times X | d(x,y)<r\}.$$ Then the collection $\{B_d(r):r>0\}$ forms a $\cU$-equivalence class on $X$ that is called the $\cU$-equivalence class generated by $d$ and denoted by $\cU_d$.

{\bf Theorem 4.6.} If a topological space can be embedded into a product of transitive pseudo-metric spaces, then it is equivalently uniformisable.

{\bf Proof:} Suppose $\{(X_i,d_i): i\in I\}$ is an indexed collection of transitive pseudo metric spaces and let $f:(X,\cT)\lora(\Pi X_i,\Pi \cT_{d_i})$ be a topological embedding.
For each $i$, let $\cU_i=\cU_{d_i}$ and $\cU=\Pi\cU_i$. Hence $\cT_{_{\cU}}=\cT_{_{\Pi\cU_i}}=\Pi\cT_{_{d_i}}$ and so
$f:(X,\cT)\lora(\Pi X_i,\cT_{_{\cU}})$ is a topological embedding. Now the result follows from proposition 4.4.$\rul$

{\bf Definition 4.7.} Let $\cD$ be a family of transitive pseudo metrics on set $X$ and let $\cS_\cD$ be the collection $\{B_d(r):r>0, d\in\cD\}$.
The $\cU$-equivalence class generated by $\cS_\cD$ is called the $\cU$-equivalence class generated by $\cD$ and denoted by  $\cU_\cD$.
Note that since for each $d\in\cD$ and each $r>0$, $d$ is $r$-transitive, then $B_d(r)$ is an equivalence relation on $X$.

{\bf Proposition 4.8.} Let $\cD$ be a collection of transitive pseudo metrics on a set $X$. Then $\cT_{_{\cD}}=\cT_{_{\cU_\cD}}$ where $\cT_{_{\cD}}$ is the topology generated by sub-base
\[\tilde{\cS_\cD}=\{B_d(x,r):x\in X, r>0, d\in\cD\}.\]

{\bf Proof:} Let $G\in \cT_{_{\cU_\cD}}$ and let $x\in G$. Then $U[x]\sse G$ for some $U\in \cU_\cD$. By the definition of $\cU_\cD$, $U=\ds\bigcap^n_{i=1} B_{d_i}(r_i)$ where $r_i>0$ and $d_i\in\cD$.  For each $i$, let $x_i=x$. Then \linebreak $x\in\ds\bigcap^n_{i=1} B_{d_i}(x_i,r_i)\sse G$. Hence $G\in\cT_{_\cD}$.
Conversely,  let $x\in G\in\cT_{_\cD}$. Then \linebreak $x\in\ds\bigcap^n_{i=1} B_{d_i}(x_i,r_i)\sse G$ where $x_i\in X$, $r_i>0$, $d_i\in\cD$.
If $U=\ds\bigcap^n_{i=1} B_{d_i}(r_i)$, then $U\in\cU_\cD$.
we contend $U[x]\sse G$. For let $y\in U[x]$, then $d_i(x,y)<r_i$ for all $i$.

On the other hand, $d_i(x,x_i)<r_i$ for all $i$, because $x\in\ds\bigcap^n_{i=1} B_{d_i}(x_i,r_i)$. By the transitivity of $d_i$, we have $d_i(x_i,y)<r_i$ for $1\leq i\leq n$. Thus $y\in\ds\bigcap^n_{i=1} B_{d_i}(x_i,r_i)\sse G$. Consequently $\cT_{_\cD}=\cT_{_{\cU_\cD}}$.$\rul$

{\bf Proposition 4.9.} Let $\cD$ be a collection of transitive pseudo metrics on a set $X$. For each $d\in\cD$, let $X_d$ be a copy of the set $X$ and let $Y=\Pi_{d\in\cD} X_d$. Then the evaluation function $f:(X,\cU_\cD)\lora(Y,\Pi\cU_d)$ defined by $f(x)(d)=x$ for all $d\in\cD$, is a $\cU$-embedding of $(X,\cU_\cD)$ into $(Y,\Pi\cU_d)$.

{\bf Proof:} Obviously $f$ is an injection. Let $U=(\pi_d\times\pi_d)^{-1}(U_d)$, $U_d=B_d(r)$ for some positive real $r$.
Then $(f\times f)^{-1}(U)=U_d\in\cU_\cD$. Hence $f$ is $\cU$-equivalently continuous. Finally, let $Z$ be the range of the function $f$ and let $r>0$ and $d\in\cD$ be given . Since $(f\times f)(B_d(r))=Z\times Z\cap(\pi_d\times\pi_d)^{-1}(B_d(r))$, then
$(f\times f)(B_d(r))\in \Pi\cU_d/Z$. Now let  $U=\ds\bigcap^n_{i=1} B_{d_i}(r_i)\in \cU_\cD$. Then $(f\times f)(U)=\ds\bigcap^n_{i=1} (f\times f)(B_{d_i}(r_i))$ because $f$ is a bijection. Hence $(f\times f)(U)\in\Pi \cU_d/Z$. Consequently, $f$ is a $\cU$-embedding.$\rul$

{\bf Theorem 4.10.} An equivalently uniformisable topological space which its associated $\cU$-equivalence class generated by a collection of transitive pseudo metrics, can be topologically embedded into a product of transitive pseudo-metric spaces.

{\bf Proof:} Let $(X,\cT)$ be the  topological space induced by $\cU_\cD$ where ${\cD}$ is a collection of transitive pseudo metrics on $X$. By Proposition 4.9, there exists a $\cU$-embedding $f:(X,\cU_\cD)\lora (\Pi X_d,\Pi\cU_d)$. Hence $f:(X,\cT_{_{\cU_\cD}})\lora (\Pi X_d,\cT_{_{\pi\cU_d}})$ is a topological embedding. But $\cT=\cT_{_{\cU_\cD}}$ and $\cT_{_{\Pi\cU_d}}=\Pi\cT_{_d}$. Thus
$f:(X,\cT)\lora (\Pi X_d,\Pi\cT_{_d})$ is a topological embedding.$\rul$

\end{document}